\documentclass[11pt]{article}

\usepackage{amsmath,amssymb,amsthm}
\usepackage{a4wide}
\usepackage[cp866]{inputenc}
\usepackage{url}
\usepackage{srcltx}

\allowdisplaybreaks

\let\tilde\widetilde
\let\epsilon\varepsilon
\let\phi\varphi
\let\upn=\textup
\def\P{\mathop{\mbox{\bf{P}}}\nolimits}
\def\nfrac#1#2{\mbox{\small$\dfrac{#1}{#2}$}}
\def\nbinom#1#2{\mbox{\small$\dbinom{#1}{#2}$}}
\renewcommand{\int}{\intop\limits}

\newtheorem*{singleremark}{Remark}
\newtheorem{conjecture}{Conjecture}

\begin{document}

\title{Minimum Number of Edges in a Hypergraph Guaranteeing\\ a Perfect Fractional
Matching and the MMS Conjecture}


\author{Vladimir Blinovsky\thanks{Supported by the S\~ao Paulo Research
Foundation (FAPESP), project nos.~2012/13341-8 and 2013/07699-0, and NUMEC/USP,
Project MaCLinC-USP.}}

\date{\small Instituto de Matem\'atica e Estat\'{\i}stica, Universidade de S\~ao
Paulo,\\ Rua do Mat\~ao 1010, 05508-090, S\~ao Paulo, Brazil\\ Kharkevich Institute for
Information Transmission Problems,\\ Russian Academy of Sciences, B. Karetnyi 19,
Moscow, Russia\\ {\tt vblinovs@yandex.ru}}

\maketitle\bigskip

\begin{abstract}
In this paper we prove the Ahlswede--Khachatrian conjecture~\cite{1} up to a finite
number of cases, which can be checked using modern computers. This conjecture implies
the conjecture from~\cite{2} and the Manickam--Mikl\'os--Singhi conjecture.
\end{abstract}

\section{Introduction}

For a ground set $[n]=\{1,\ldots,n\}$, a (uniform) \emph{hypergraph\/} $H(n,k)=(n,E)$
is $[n]$ together with a subset of edges $E\subset \nbinom{[n]}{k}$. A \emph{perfect
fractional matching\/} of $H$ is a set of nonnegative real numbers
$(\alpha_1,\ldots,\alpha_{|E|})$, $\alpha_j\ge 0$, $\sum\limits_{j=1}^{|E|}\alpha_j
=1$, $\sum\limits_{e\in E}\alpha_e e=\Bigl(\nfrac{k}{n},\ldots,\nfrac{k}{n}\Bigr)$. A
hypergraph need not have a perfect fractional matching. In~\cite{1} Ahlswede and
Khachatrian proposed the following conjecture.

\begin{conjecture}\label{co1}
Let\/ $\mathcal{H}_f$ be the set of hypergraphs that do not have perfect fractional
matchings. Then
\begin{equation}\label{e1}
p(n,k)\triangleq\max_{H\in\mathcal{H}_f}|E| =\max_{k\ge s\ge
1}\sum_{i=0}^{k-s}\binom{n_s}{i+s}\binom{n-n_s}{k-s-i},
\end{equation}
where $n_s =\lceil ns/k\rceil -1$.
\end{conjecture}

In~\cite{1} this conjecture was given in other, less standard terms (cone
dependence), and this drew less attention of specialists to this problem.

We relax the conditions of this conjecture a little and rewrite it as
\begin{equation}\label{e2}
p(n,k)=\max_{n-1\ge a\ge 1}\sum_{i>ka/n}\binom{a}{i}\binom{n-a}{k-i}.
\end{equation}
To see that~\eqref{e1} is equivalent to~\eqref{e2}, one can observe that the above
sum decreases when $a$ decreases from $n_s$ to $n_{s-1}-1$.

Now we turn to another problem. Let
$$
(\beta_1,\ldots,\beta_n),\quad \beta_j\in\mathbb{R}^1,\quad \sum_{j=1}^{n}\beta_j =0,
$$
and
$$
U(\{\beta_j\})=\biggl\{e\in\binom{[n]}{k}:\: \sum_{j\in e}\beta_j\ge 0\biggr\}.
$$
Define
$$
q(n,k)=\min_{\{\beta_j\}}|U(\{\beta_j\})|.
$$
The following conjecture was proposed in~\cite{2}.

\begin{conjecture}\label{co2}
The relation
\begin{equation}\label{e3}
q(n,k)=\min_{n-1\ge a\ge 1}\sum_{i\ge ka/n}\binom{a}{i}\binom{n-a}{k-i}
\end{equation}
is valid.
\end{conjecture}

A close relation between the problem of determining $p(n,k)$ and that of determining
$q(n,k)$ was also shown in~\cite{2}. In particular, it was shown that Conjecture~1
implies Conjecture~2. Actually, these two problems are equivalent. We will show this
later.

Finally, we come to the Manickam--Mikl\'os--Singhi (MMS) conjecture \cite{3,4,5,6}:
\begin{equation}\label{e4}
q(n,k)=\binom{n-1}{k-1},\quad n\ge 4k.
\end{equation}
Thus, to prove the MMS conjecture (assuming that Conjecture~2 is true), we have to
show that if $n\ge 4k$, then
\begin{equation}\label{e88}
\min_{n-1\ge a\ge 1}\sum_{i\ge ka/n}\binom{a}{i}\binom{n-a}{k-i} =\binom{n-1}{k-1}.
\end{equation}
This will be done in the Appendix.

Taking into account a natural bijection between subsets of $[n]$ and binary
$n$-tuples, we will identify them in what follows.

Now we prove that both problems are equivalent. This will imply that the proof of
Conjecture~1 follows from the proof of Conjecture~2.

Consider $A\subset\nbinom{[n]}{k}$ as a set of vertices of the hypersimplex $\Gamma
(n,k)\subset \mathbb{R}^n$ such that its convex hull $K(A)$ does not contain the
point $\smash[t]{\Bigl(\nfrac{k}{n},\ldots,\nfrac{k}{n}\Bigr)}$. Then there exists a
hyperplane
$$
\sum_{j=1}^{n}\omega_j y_j =0,
$$
such that $\sum\limits_{j=1}^{n} \omega_j =0$ (it contains the center
$(k/n,\ldots,k/n)$) and $K(A)$ belongs to one of the open half-spaces into which the
hyperplane separates $\mathbb{R}^n$. Vice versa, if there exists a hyperplane with
the above properties, then $K(A)\not\ni (k/n,\ldots,k/n)$. If $A$ is maximal, then
$$
M=\Biggl\{x\in\binom{[n]}{k}:\: \sum_{j=1}^{n}\omega_j x_j\ge 0 \Biggr\}
$$
is minimal.

This shows the equivalence of the problems.

Assume that $\omega_1\ge\omega_2\ge\ldots\ge\omega_{n}$,
$\sum\limits_{j=1}^{n}\omega_j =0$. Then the space of such $\omega$ has a basis
$\{z_j,\: j=1,\ldots,n-1\}$ of the form $z_j=(n-j,\ldots,n-j,-j,\ldots,-j)$. Each $y$
from this space has a representation
$$
y=\sum_{j=1}^{n-1}y_j z_j
$$
with nonnegative coefficients $y_j\ge 0$. Fix $x\in\nbinom{[n]}{k}$. Then (letting
$y_i =\sum\limits_{j=1}^{n-1}z_{ji}y_j$) we have
$$
\begin{aligned}
(x,y)&=\sum_{j=1}^{n}x_i y_i =\sum_{j=1}^{n}x_i \sum_{j=1}^{n-1}z_{ji}y_j=
\sum_{j=1}^{n-1}y_j \sum_{i=1}^{n}x_i z_{ji}=\sum_{j=1}^{n-1}y_j
\Biggl(n\sum_{i=1}^{j}x_i -jk\Biggr)\\ &= n\sum_{j=1}^{n-1}y_j \sum_{i=1}^{j}x_i
-k\sum_{j=1}^{n-1}jy_j.
\end{aligned}
$$
Again using the fact that the above conditions are homogeneous and dividing the
right-hand side of the last chain of equalities by $\sum\limits_{j=1}^{n-1}jy_j$, we
obtain the condition
$$
\sum_{j=1}^{n-1}\frac{y_j}{\sum\limits_{j=1}^{n-1}j y_j}\sum_{i=1}^{j}x_i
\ge\frac{k}{n}.
$$
The last inequality is equivalent to
$$
\sum_{j=1}^{n-1}\gamma_j x_j\ge\frac{k}{n},\quad \gamma_j\ge 0,\quad
\sum_{j=1}^{n-1}\gamma_j =1.
$$
If we consider only the strict inequality
$$
\sum_{j=1}^{n-1}\gamma_j x_j >\frac{k}{n},\quad \gamma_j\ge 0,\quad
\sum_{j=1}^{n-1}\gamma_j =1,
$$
and try to find the maximum number of its solutions for
$x\in\smash[b]{\nbinom{[n]}{k}}$, then this is exactly the problem of determining
$p(n,k)$.

Hence it follows that
\begin{equation}\label{e6}
p(n,k)=\max_{\{\gamma_j\}}{}\Biggl| x\in\binom{[n]}{k}:\: \sum_{j=1}^{n-1}x_j
\gamma_j >\frac{k}{n}\Biggr|.
\end{equation}

\paragraph{History.}

As was already noted, the problem of determining $p(n,k)$ was first posed by Ahlswede
and Khachatrian~\cite{1}. From their results, using considerations from~\cite{2}, it
can easily be deduced that
$$
q(n,k)=\binom{n-1}{k-1},\quad n\ge 2k^3.
$$
In~\cite{8} this equality was proved for $n\ge\min\{2k^3,33k^2\}$. Finally, validity
of this equality for $n\ge10^{46}k$ was proved in~\cite{9}.

\section{Proof of Conjecture~1}

Now assume that $k<n$. We use equality~\eqref{e6}. Consider the following
function:
$$
f(\{\gamma_1,\ldots,\gamma_ {n-1}\})=\frac{1}{\sqrt{2\pi}}\sum_{x\in\binom{[n]}{k}}
\int_{-\infty}^{\bigl(\sum\limits_{j=1}^{n-1}\gamma_j x_j
-\frac{k}{n}\bigr)/\sigma}e^{-\frac{z^2}{2}}\,dz.
$$
Define
$$
N(\gamma_1,\ldots,\gamma_{n-1})=\Biggl| x\in\binom{[n]}{k}:\:
\sum_{j=1}^{n-1}\gamma_j x_j >\frac{k}{n}\Biggr|.
$$
Then
$$
|N(\{\gamma_j\})-f(\{\gamma_j\})|<\epsilon (\sigma),\quad \epsilon
(\sigma)\xrightarrow{\sigma\to 0\,} 0,
$$
uniformly over $\{\gamma_j\}$ such that
\begin{equation}\label{e44}
\Biggl| \sum_{j=1}^{n-1}\gamma_j x_j -\frac{k}{n}\Biggr| >\delta,\quad \forall
x\in\binom{[n]}{k}.
\end{equation}
For extremal $\gamma$, i.e., those with $N(\gamma)=p(n,k)$, it is easy to see that
$\gamma$ satisfies condition~\eqref{e44} for some $\delta >0$, because otherwise, if
$\sum\limits_{j=1}^{n-1}\tilde{\gamma}_j x^0_j =\nfrac{k}{n}$ for some
$x^0\in\nbinom{[n]}{k} $, there exists $\gamma'$ close to $\tilde{\gamma}$ that does
not violate the conditions
$$
\gamma_j'\ge 0,\quad \sum_{j=1}^{n-1} \gamma_j' =1,\quad \sum_{j=1}^{n-1}\gamma'_j
x_j >\frac{k}{n}
$$
and $\sum\limits_{j=1}^{n-1}\gamma'_j x^0_j >\nfrac{k}{n}$.

Hence, assuming that we are interested in extremal $\gamma$, we may assume
that~\eqref{e44} is satisfied.

Now assume without loss of generality that $\gamma_1\ge\ldots\ge\gamma_{n-1}$. Since
we have the constraints $\gamma_j\ge 0$, we should look for the extremum among
$\gamma$ such that
$$
\gamma_{a+1}=\ldots= \gamma_{n-1}=0,\quad a =1,\ldots, n-1
$$
($a=n-1$ means that we are not imposing any zero condition on $\gamma$). Assume that
this condition is valid for some $a\ge 5$. (The cases with $a\le 4$ are easy to
consider.) Then, since $\gamma_a=1-\sum\limits_{j=1}^{a-1}\gamma_j$, we have
\begin{equation}\label{e45}
f'_{\gamma_j}=\frac{1}{\sqrt{2\pi}\sigma}\sum_{x\in\binom{[n]}{k}:\: j\in x,\:
a\notin x}e^{-\bigl(\sum\limits_{j=1}^{a}\gamma_j x_j
-\frac{k}{n}\bigr)^2/(2\sigma^2)} -\frac{1}{\sqrt{2\pi}\sigma}
\sum_{x\in\binom{[n]}{k}:\: j\notin x,\: a\in x}
e^{-\bigl(\sum\limits_{j=1}^{a}\gamma_j x_j -\frac{k}{n}\bigr)^2/(2\sigma^2)}.
\end{equation}

Next we show that we may assume that these equalities can be jointly attained on step
functions $\gamma_j=\gamma_a$ with $j\in [a]$. Indeed, choose a parameter $\sigma$
sufficiently small and then fix it. Then, to satisfy equations~\eqref{e45}, we should
assume the equalities
$$
\sum_{x\in\binom{[n]}{k}:\: j\in x,\: a \notin x}e^{-{\textstyle(}(\gamma
,x)-\frac{k}{n} {\textstyle)}^2/(2\sigma^2)} =\sum_{x\in\binom{[n]}{k}:\: a\in x,\: j
\notin x}e^{-{\textstyle(}(\gamma,x)-\frac{k}{n} {\textstyle)}^2/(2\sigma^2)}.
$$
For these equalities to be satisfied, we should assume that the exponents from the
left-hand sum are equal to the corresponding exponents from the right-hand sum; i.e.,
for each given $j\in [a -1]$,
\begin{equation}\label{er89}
\Bigl((\gamma,x)-\frac{k}{n} \Bigr)^2 =\Bigl((\gamma,y)-\frac{k}{n} \Bigr)^2,
\end{equation}
where $x\in\nbinom{[n]}{k}$, $j\in x$, $y\in\nbinom{[n]}{k}$, $a\in y$, and where
$x\setminus \{j\}$ and $y\setminus \{a\}$ run over all sets of cardinality $k-1$ from
$[n-j-a]$. We rewrite conditions~\eqref{er89} as follows:
\begin{multline*}
\gamma_j^2 +(\gamma_{j_1}+\ldots+\gamma_{j_{k-1}})^2 -2\frac{k}{n} \gamma_j
-2\frac{k}{n} (\gamma_{j_1}+\ldots+\gamma_{j_{k-1}})
+\gamma_j(\gamma_{j_1}+\ldots+\gamma_{j_{k-1}})\\ =\gamma_a^2
+(\gamma_{m_1}+\ldots+\gamma_{m_{k-1}})^2 -2\frac{k}{n} \gamma_a -2\frac{k}{n}
(\gamma_{m_1}+\ldots+\gamma_{m_{k-1}})+\gamma_a(\gamma_{m_1}+\ldots
+\gamma_{,m_{k-1}}).
\end{multline*}
Summing both sides of this equality over all admissible choices of
$j_1,\ldots,j_{k-1}$ and $m_1,\ldots,m_{k-1}$, we obtain
\begin{equation}\label{er56}
\binom{n-2}{k-1}\hskip-1pt\biggl(\!\gamma_j^2 -2\frac{k}{n} \gamma_j \!\biggr)\!
-2\frac{k}{n} R +2\gamma_jR =\binom{n-2}{k-1}\hskip-1pt\biggl(\!\gamma_a^2
-2\frac{k}{n} \gamma_a\!\biggr)\! -2\frac{k}{n} R +2\gamma_aR,
\end{equation}
where
$$
R= \sum_{x\in\binom{[n]\setminus\{j, a\}}{k-1}} (\gamma,x)=\binom{n-3}{k-2}\sum_{m\ne
j,a}\gamma_{m}=\binom{n-3}{k-2}(1-\gamma_j-\gamma_a).
$$
Denote by $\mu$ the term $\nfrac{k}{n}$ in~\eqref{er56}. From~\eqref{er56} it follows
that $\gamma_j$ can take at most two values:
\begin{equation}\label{er91}
\begin{gathered}
\gamma_j=\gamma_a,\\ \gamma_j+\gamma_a=\lambda \triangleq
2\frac{\mu-\nfrac{k-1}{n-2}}{1-2\nfrac{k-1}{n-2}}.
\end{gathered}
\end{equation}
Note that we can vary $\mu$ a little without violating all the above considerations.
Next we show how we can eliminate the possibility that $\gamma_j$ takes the second
value. Assume first that to each $x$ such that $ |x\cap [a]|=p$ there corresponds
some $y$ such that $|y\cap [a]|=p$ for all $x\in\nbinom{[n]}{k}$ and $p$. For a given
$p$, we sum the left- and right-hand sides of~\eqref{er89} over $x$ and the
corresponding $y$ such that $|x\cap [a]|=p$. Then, similarly to the case of summation
over all $x$, we obtain two possibilities: either
$$
\gamma_j=\gamma_a
$$
or
\begin{equation}\label{et1}
\gamma_j+\gamma_a=2\frac{\mu -\nfrac{p-1}{a-2}}{1-2\nfrac{p-1}{a-2}}.
\end{equation}
Since we can vary $p$, it follows that equality~\eqref{et1} for some $p$ would
contradict the second equality from~\eqref{er91}.

Now assume that for some $b$ we have
\begin{equation}\label{er77}
\gamma_j=
\begin{cases}
\lambda-\gamma_a,& j\le b,\\ \gamma_a,& j\in[b+1,a].
\end{cases}
\end{equation}
Since $\smash[b]{\sum\limits_j \gamma_j}=1$, we have the following condition on
$\gamma_a$ and $\nfrac{k}{n}$:
\begin{equation}\label{ek1}
b\lambda +(a -2b)\gamma_a =1.
\end{equation}
Let $\gamma_j=\lambda -\gamma_a$. Assume also that to some $x$ such that $|x\cap
[a]|=p$ there corresponds some $y$ such that $|y\cap [a]|=q$ for some $p\ne q$.
From~\eqref{er89} it follows that there exist two possibilities: either
$$
(\gamma,x)=(\gamma,y)
$$
or
\begin{equation}\label{e34}
(\gamma,x)+(\gamma,y)=2\frac{k}{n}.
\end{equation}
Each of these equalities imposes some condition; the first equality, the condition
(for some integers $p_1$ and $p_2$)
$$
p_1 \gamma_a+ p_2 \lambda =0,
$$
which is either inconsistent with equality~\eqref{ek1} or together with
equality~\eqref{ek1} uniquely determines the value of $\mu$.

On the other hand, equality~\eqref{e34} imposes the condition (for some integers
$p_3$ and $p_4$)
\begin{equation}\label{ed1}
p_3 \gamma_a+p_4 \lambda =2\frac{k}{n}.
\end{equation}
It is possible that equality~\eqref{ek1} together with equality~\eqref{ed1} does not
determine the value $k/n$. In this case there are the following two possibilities.
The first is that there exist $x: |x\cap [a]|=m$ (where $m$ can be equal to $p$ or
$q$) and the corresponding $y:\: |y\cap [a]|=v$ with $v\ne p,q$.

The second possibility is that each $x: |x\cap [a]|=m$ (here $m\ne p,q$) corresponds
to $y: |y\cap [a]|=m$. In this second case we return to the case which leads to
equalities~\eqref{et1} (because if $a\ge 5$, the number of such $m\ne p,q$ is greater
that $1$).

If we have the first possibility, then we have one additional equation
\begin{equation}\label{el1}
q_3 \gamma_a+q_4 \lambda =2\frac{k}{n},
\end{equation}
which is either inconsistent with equalities~\eqref{ek1} and~\eqref{ed1} or together
with them determines a unique value of $\nfrac{k}{n}$.

We see that if $b>1$ and $\gamma_j=\lambda -\gamma_a$ with $j\le b$, then $k/n$ can
take values from some discrete finite set only. Considering other numbers $\mu$
instead of $k/n$ which are sufficiently close to $k/n$, we can achieve a situation
where none of the values from this finite set is equal to $\mu$. Again we emphasize
that such variation can always be made without violating relation~\eqref{e44}.

Let $N(\gamma)$ achieve its extremum at $\bar{\gamma}$, and $f(\gamma)$, at
$\tilde{\gamma}$. We have
\begin{gather*}
|N(\tilde{\gamma})-f(\tilde{\gamma})|<\epsilon,\\
|N(\bar{\gamma})-f(\bar{\gamma})|<\epsilon.
\end{gather*}
Then
$$
N(\bar{\gamma})<f(\bar{\gamma})+\epsilon<f(\tilde{\gamma})+\epsilon
<N(\tilde{\gamma})+2\epsilon.
$$
Now $N(\gamma)$ is a positive integer, and the last inequalities mean that
$$
N(\tilde{\gamma})=N(\bar{\gamma}).
$$
Hence follows Conjecture~1.

\appendix{}

Fist, note the fact that if $k\mid n$, then $q(n,k)=\nbinom{n-1}{k-1}$. This
easily follows from the lemma in~\cite{7}.

Next, we will use the fact that if
$$
p(n,k)=\binom{n-1}{k},
$$
then
$$
p(n+k,k)=\binom{n+k-1}{k}.
$$
A proof using double counting argument can be found in~\cite{10} (where only $q(n,k)$
was considered, but the problem for $p(n,k)$ is equivalent). Thus, we may assume that
$k\in (n/5,n/4)$.

Next we assume that $n\nmid ka$. We estimate the probability $\P(k\le ka/n)$ in
several steps. First, we use the Berry--Esseen inequality for a hypergeometric
distribution. We will use considerations from~\cite{12}. The problem is that
in~\cite{12} there was not computed the constant $C$ in the inequality
$$
\left| \P\left(a<\frac{i-ka/n}{\sigma}\le b\right) -(\Phi (b)-\Phi (a))\right|
<\frac{C}{\sigma},
$$
where
\begin{equation}\label{e00}
\sigma^2 =\frac{ka}{n}\left(1-\frac{a}{n}\right)\Bigl(1-\frac{k}{n}\Bigr).
\end{equation}
Thus, we have to repeat those considerations in a way that allows us to obtain an
appropriate upper bound for $C$. For details, we refer the reader to~\cite{12}; we
are only interested in estimates for the constant in Theorem~2.2.

Let $\delta =1/20$ and $\sigma >55$, $n>12\cdot 10^4$,
$\phi(x)=\nfrac{1}{\sqrt{2\pi}}\exp (-x^2 /2)$,
$\Phi(x)=\nfrac{1}{\sqrt{2\pi}}\int_{-\infty}^x e^{-z^2 /2}\,dz$,
$$
\begin{aligned}
K_2 &=\max\{i\in\mathbb{Z}_+:\: \tilde{x}(i)\ge -\delta\sigma\},\\ K_1 &=
\min\{i\in\mathbb{Z}_+: \tilde{x}(i)\ge -1\},\\ K_0& =\max\{i\in\mathbb{Z}_{+}:\:
\tilde{x}(i)\le 0\},
\end{aligned}
$$
where
$$
\tilde{x}(i)=\frac{i-ka/n}{\sigma}.
$$
Then \cite[equation~(3.15)]{12}
$$
\Delta \triangleq \P(i<K)+\sum_{i=K}^{\lfloor ka/n\rfloor} \left|
P(i)-\frac{1}{\sigma}\phi (\tilde{x} (i))\right| +\Biggl|\sum_{j=K}^{\lfloor
ka/n\rfloor}\frac{1}{\sigma}\phi (\tilde{x}(j))-\Phi (x)\Biggr|=I_1 +I_2 +I_3.
$$
Since $K_2 -1<ka/n -\delta\sigma^2\le K_2$, by using the Chebyshev inequality we
obtain
$$
\begin{aligned}
I_1 &=\P(i\le K_2 -1)\\ &\le\P\left(\left|\frac{i-ka/n}{\sigma}\right|
\ge\left|\frac{K_2-ka/n-1}{\sigma}\right|\right)\\ &\le\frac{\mathop{\rm
Var}(i)}{(K_2 -1-ka/n)^2}\le\frac{n\sigma^2}{n-1}(\delta\sigma^2)^{-2}<0.1323.
\end{aligned}
$$
Taking into account that $k\in (n/5,n/4)$, we have (see \cite[inequality~(3.20)]{12})
\begin{multline*}
\sum_{i=K_2}^{\lfloor ka/n\rfloor}\left| P(i)-\frac{1}{\sigma}\phi
(\tilde{x}(i))\right|\\[-3pt]
\begin{aligned}
&\le\frac{17}{3\sqrt{2\pi}\sigma^2}\exp\{\sigma^{-2}\}\Biggl(2\int_0^\infty x^3 \exp
(-0.07 x^2)\,dx +\left(\frac{4(3/(2\cdot 0.07))^{3/2}}{\sigma}\right)\exp
(-3/2)\Biggr)\\ & <\frac{17}{3\sqrt{2\pi}\sigma^2}\exp\{\sigma^{-2}\}\cdot
16<\frac{60}{\sigma^2}\exp\{\sigma^{-2}\}.
\end{aligned}
\end{multline*}
This inequality differs from~(3.20) in \cite{12}, because we transform it with the
use of the inequality
$$
\sum_{i=-\infty}^{\infty} g(ih)=2\sum_{i=0}^{\infty}g(ih)\le
2\int_0^{\infty}g(x)\,dx+4hg(x_0),
$$
which is true for a symmetric nonnegative unimodal function $g$ on $[0,\infty)$ with
a maximum (one of two) at $x_0$ and for $h>0$. We use it to approximate the sum
in~(3.20) as follows:
$$
\sum_{i=K_2}^{\lfloor ka/n\rfloor} |\tilde{x}(i)|^3 \exp (-0.07\tilde{x}^2(i))
\le\int_{-\infty}^{\infty}|x|^3\exp (-0.07x^2) \,dx+\frac{4(3/(2\cdot
0.07))^{3/2}}{\sigma}\exp (-3/2)< 16.
$$

Furthermore (see \cite[relations~(3.21)]{12}), we have
$$
\sum_{i=K_1}^{K_0} \left| P(i)-\frac{1}{\sigma}\phi
(\tilde{x}(i))\right|\le\frac{85}{12\sqrt{2\pi}\sigma}\exp\{\sigma^{-1}\}
<\frac{3}{\sigma}\exp\{\sigma^{-1}\}.
$$
Thus,
$$
I_2\le\frac{60}{\sigma^2}\exp\{\sigma^{-2}\}
+\frac{3}{\sigma}\exp\{\sigma^{-1}\}<0.077.
$$
Finally, for $I_3$ we have the estimate \cite[ineqiuality~(3.22)]{12}
$$
\begin{aligned}
I_3 &\le\frac{1}{12\sigma^2}\left[\frac{1}{\sqrt{\pi}}+1 +\frac{10}{\sqrt{2\pi}}
\exp\{-1/(8\sigma^2)\right]+\Phi (1/(2\sigma)-\Phi (-1/(2\sigma))+\Phi (-\delta\sigma
+1/(2\sigma))\\ &\le\frac{1}{12\sigma^2}
\left[\frac{1}{\sqrt{\pi}}+1+\frac{10}{\sqrt{2\pi}} \exp\{-1/(8\sigma^2)\right]
+\frac{1}{\sqrt{2\pi}\sigma}+\frac{\exp (-(\delta\sigma -1/(2\sigma))^2
/2)}{\sqrt{2\pi}(\delta\sigma -1/(2\sigma))}\\ &<0.011.
\end{aligned}
$$
Provided that $\sigma>55$ and $n>12\cdot 10^4$, from these considerations we finally
obtain
$$
I_1 +I_2 +I_3 <0.2203.
$$
Hence, the Berry--Esseen inequality is of the form
$$
\left| \P(i\le ka/n)-\frac{1}{2}\right| =\left|\P(i>ka/n)
-\frac{1}{2}\right|<0.2203<1/4.
$$
Thus, in this case it follows that for $k\in (n/5,n/4)$ we have
\begin{equation}\label{et}
\P(i>ka/n)\binom{n}{k}\le\binom{n-1}{k}.
\end{equation}
Using formula~\eqref{e00} for $\sigma$, one can easily see that the condition
$\sigma\ge 55$ is satisfied if $a\in [n/5,n-n/5]$ and $n>12\cdot 10^{4}$. The next
consideration allows us to reduce possible values of $a$ for which~\eqref{et} is
valid to $a\in [C,n-C]$, for some constant $C$. Using Stirling's formula, it is easy
to see that for $i\le ka/n$ we have the inequality
$$
\frac{\nbinom{n-a}{k-i}}{\nbinom{n}{k}}
<2^{(n-a)H\left(\frac{k-i}{n-a}\right)-nH\left(\frac{k}{n}\right)}
\frac{\sqrt{\nfrac{k}{n}\left(1-\nfrac{k}{n}\right)}}
{\sqrt{\left(\nfrac{k-i}{n-a}\right)\left(1-\nfrac{k-i}{n-a}\right)
\left(1-\nfrac{a}{n}\right)}} e^{\frac{1}{12n}
+\frac{1}{12(k-i)}+\frac{1}{12(n-a-k+i)}}.
$$

Since we may assume that $i<ka/n$ with $n>12\cdot 10^4$ and $a\le n/5$, we obtain the
estimate
$$
\frac{\nbinom{n-a}{k-i}}{\nbinom{n}{k}}\le\left(\frac{k}{n}\right)^i
\left(1-\frac{k}{n}\right)^{a-i} \left(1-\frac{a}{n}\right)^{-1/2}e^{10^{-3}}\le
1.1203 \left(\frac{k}{n}\right)^i \left(1-\frac{k}{n}\right)^{a-i}.
$$
Here we used the fact that for $i\le ka/n$ we have
$$
(n-a)H\left(\frac{k-i}{n-a}\right)-H\left(\frac{k}{n}\right)\le
i\ln\frac{k}{n}+(a-i)\ln\left(1-\frac{k}{n}\right).
$$
Hence,
\begin{equation}\label{ep}
\sum_{i<ka/n}\binom{a}{i}\binom{n-a}{k-i}\le
1.1203\sum_{i<ka/n}\binom{a}{i}\left(\frac{k}{n}\right)^i
\left(1-\frac{k}{n}\right)^{a-i} \binom{n}{k}.
\end{equation}
To estimate the sum
$$
\sum_{i<ka/n}\binom{a}{i}\left(\frac{k}{n}\right)^i \left(1-\frac{k}{n}\right)^{a-i}
$$
we note that this is the probability $P_0$ that the sum of $a$ i.i.d.\ Bernulli
variables exceeds the average and use the Berry--Esseen inequality~\cite{11}
(relaxing the coefficients a little):
$$
\left| P_0 -\frac{1}{2}\right| <\frac{\rho +0.43\sigma_1^3}{3\sigma_1^3 \sqrt{a}},
$$
where $\sigma_1^2 =\nfrac{k}{n}\left(1-\nfrac{k}{n}\right)\ge\nfrac{5}{25}$ and $\rho
=\sigma_1^2 (1-2\sigma_1^2)$. Thus,
$$
\left| P_0 -\frac{1}{2}\right| <\frac{0.71}{\sqrt{a}}
$$
ans
$$
\begin{aligned}
\sum_{i<ka/n}\binom{a}{i}\binom{n-a}{k-i}&\le1.026
\sum_{i>ka/n}\binom{a}{i}\left(\frac{k}{n}\right)^i \left(1-\frac{k}{n}\right)^{a-i}
\binom{n}{k}\\ &<1.1203\left(\frac{1}{2} +\frac{0.71}{\sqrt{a}}\right)\binom{n}{k}.
\end{aligned}
$$
The right-hand side of the last inequality is less than $\nbinom{n-1}{k}$ if
$$
1.1203\left(\frac{1}{2} +\frac{0.71}{\sqrt{a}}\right) <\frac{3}{4},
$$
i.e., if $a>14$. Next, if $a=n-b>n-n/5$, we can rewrite the sum on the left-hand side
of~\eqref{ep} as
$$
\sum_{i<kb/n}\binom{b}{i}\binom{n-b}{k-i}.
$$
To estimate the ratio
$
{\nbinom{n-b}{i}}\mathrel{\Bigl/}{\nbinom{n}{k}},
$
repeat the above arguments for $b$ instead of $a$ to obtain
$$
\sum_{i<kb/n}\binom{b}{i}\binom{n-b}{k-i}\le 3/4\binom{n}{k}
$$
if $b=n-a<9$.

\begin{singleremark}
For $k\le n/4$ and\/ $n<12\cdot 10^4$\upn, the inequality
$$
\sum_{i>ka/n}\binom{a}{i}\binom{n-a}{k-i}\le\binom{n-1}{k}
$$
can be checked using the Wolfram\/ \emph{Mathematica\/}$^\circledR$ software\upn, but
this requires a fast computer. For $a\le 14$\upn, $a\ge n-14$\upn, $k\le n/4$\upn,
this can be done by hand.
\end{singleremark}

\end{document}